# An Improved NSGA-II with local search for multi-objective energy-efficient flowshop scheduling problem

Vigneshwar Pesaru[1,2*], Venkataramanaiah Saddikuti[3]

[1]Department of Manufacturing and Systems Engineering, University of Texas, Arlington at Arlington, TX 76019, United States.
[2]Bank of America-Enterprise Operations Research, 7105 Corporate Dr, Plano, TX 75024, United States.
[3]Department of Operations Management, Indian Institute of Management-Lucknow, NOIDA Campus, Uttar Pradesh 201307, India.

**Abstract:** There has been an increasing concern to reduce the energy consumption in manufacturing and other industries. Energy consumption in manufacturing industries is directly related to efficient schedules. The contribution of this paper includes: i) a permutation flowshop scheduling problem (PFLSP) mathematical model by considering energy consumed by each machine in the system. ii) an improved non-dominated sorted genetic algorithm with Taguchi method with further incorporating local search (NSGA-II_LS) is proposed for the multi-objective PFLSP model. iii) solved 90 benchmarks problems of Taillard (1993) for the minimisation of flowtime (FT) and energy consumption (EC). The performance of the proposed NSGA_LS algorithm is evaluated on the benchmark problems selected from the published literature Li et. al, (2018). From these results, it is noted that the proposed algorithm performed better on both the objectives i.e., FT and EC minimization in 5 out of 9 cases. On FT objective our algorithm performed better in 8 out of 9 cases and on EC objective 5 out of 9 cases. Overall, the proposed algorithm achieved 47% and 15.44% average improvement in FT and EC minimization respectively on the benchmark problems. From the results of 90 benchmark problems, it is observed that average difference in FT and EC between two solutions is decreasing as the problem size increases from 5 machines to 10 machines with an exception in one case. Further, it is observed that the performance of the proposed algorithm is better as the problem size increases in both jobs and machines. These results can act as standard solutions for further research.

## 1. Introduction

Extensive research has been conducted on single and multi-objective optimization of scheduling problems in a general flowshop or job shop with the aim to minimize makespan, flowtime, energy consumption, etc. With the increase in consumer desire to purchase goods and services, the demand for continuous manufacturing and improvement has become primary task. Indeed, a significant proportion of energy used in the manufacturing enterprises is mainly generated from the fossil fuels. Hence, need for minimizing the usage of fossil fuels or in other words, minimizing

*Contact: mallanjulap1@gmail.com*

the energy consumption is the important factor to be considered in manufacturing (Afshin Mansouri & Aktas, 2016). Moreover, in the recent times there has been increasing concern on the carbon efficiency of the manufacturing industry.

Since the carbon emissions in the manufacturing sector is directly related to the energy minimization. It is also considered as efficacious way to improve carbon efficiency in any industrial plant to design scheduling strategies aiming to reduce the energy minimization (Ding et al., 2016); (Nilakantan et al., 2017). Researchers developed a framework for energy efficient scheduling (EES) in order to minimize the energy based on the iterative methodology(Gahm et al., 2016). Energy efficiency has become more and more critical for the success of manufacturing companies because of rising energy prices and public perceptions of environmental conscious operations. In many cases final energy sources are not directly consumed by the production resources and thus must transform by the conversion units into applied energy sources (Rager et al., 2015).

On average, energy industries generate 28% of greenhouse gases emissions in OECD countries, followed by transport (23%), manufacturing industries (12%), agriculture (10%), industrial processes (7%) and waste (3%) (Ghazouani et al., 2021). In UK, industry electricity consumptions accounts for 31 percent of the total energy consumption. It is obliged that manufacturing companies to put more efforts in reducing their environmental impact. One way to do the operations in energy efficient manner is by shutting down the machines when they are not in use (Duflou et al., 2012). About one-half of the world's total energy consumption is contributed by the industrial sector. Thus, manufacturing enterprises have become a major source of global warming and their carbon footprints are likely to be restricted by taxes and related regulations in the future. Managing the operating costs is crucial-especially energy intensive industries such as chemical, textiles, or food(Kleindorfer et al., 2005; Fang et al., 2011a;Ngai et al., 2012;Akbar & Irohara, 2018). Because industry is acting to fulfil the growing demand for goods and consequently is one of the primary consumers of energy. In 2012, industry accounted nearly 24.2 percent of energy consumptions in the European union and it is considered to be the high time for sustainable manufacturing (Haapala et al., 2013;Mouzon & Yildirim, 2008;Yüksel et al., 2020a;S. Li et al., 2018;Saddikuti & Pesaru, 2019). By considering these major issues, we propose a model and solution approach based on NSGA-II for the flowshop scheduling problem. The contributions of this paper include the following:

1. Formulated multi-objective mathematical model by considering FT and EC in flowshop.
2. Developed multi-objective evolutionary algorithm (NSGA-II_LS) by conducting exhaustive experiments through design of experiments to identify best approximation of Pareto optimal solutions.
3. Further, evaluated the performance of the proposed NSGA-II_LS algorithm with simple NSGA-II on nine datasets from Li et al., (2018) and compared the solutions.
4. Considered 90 benchmark problems from Taillard (1993) and implemented the proposed NSGA-II_LS algorithm.
5. Further, presented results that can act as standard solutions when considered the FT and EC consumption minimization in flowshop problems when applied other heuristic or metaheuristic algorithms.

*Contact: mallanjulap1@gmail.com*

In this paper, a methodology has been proposed to optimize the scheduling of operations based on firm's requirement, i.e., by considering different aspects like energy consumption and flowtime. Most of our study focussed on the two literature lines i.e., multi-objective decision-making criteria and sustainability. The organization of the paper is as follows. Section 2 examines the appropriate literature. Section 3 constructs the mathematical model for FT and EC minimization. The constructive of NSGA-II with local search for multi-objective scenario for the permutation flowline scheduling problem is illustrated in Section 4. The experimental setup is presented in the section 5 followed by presentation and discussion of results in section 6. Section 7 concludes the work done along with managerial insights.

## 2. Summary of Literature

Incorporating the energy efficient scheduling has become a major focus in the recent times. Past research has prioritized energy optimization coupled with time efficiency within a single machine (Mouzon & Yildirim, 2008). Energy consumption and processing time of a computer control machines can vary significantly by changing cutting speed, feed rate, depth of cut, and nose radius(Ahilan et al., 2013). Machine time efficiency correlates with demand on energy. The processing speed of a machine tool is a critical factor in its energy use (Diaz et al., 2011). The relationship between the turning parameters and the power consumed based on the various machining levels and reported positive relationship between them. It is also estimated in using this relationship in scheduling too in very explicit manner (Fang et al., 2011b;G.-S. Liu et al., 2013). Construction of two MILP (mixed integer linear programming) models by taking total tardiness and energy as objective functions and solved by using C++. Computational results with randomly generated instances demonstrate that the model using assignment constraints is much more efficient than that with dichotomous constraints(Che et al., 2015). Considering variable energy prices during a day, a mathematical model is proposed to minimize energy consumption costs for single machine production scheduling during production processes (Shrouf et al., 2014). A decision-making solution presented by considering both production and energy efficiency of the unrelated parallel machine scheduling problem using the weighted sum objective of production scheduling and electricity usage(Moon et al., 2013). In fact, past review of the literature shows an integrated methodology monitoring by sensory systems, human operators and scheduling model by incorporating energy consumption paired with environmental impact consideration (Mourtzis et al., 2016).

Energy consumption and makespan time as an objective a mathematical model was developed and used as an objective function. Existing literature demonstrated the performance with respective energy utilization and time efficiency of the Ant colony optimization algorithm (Liang et al., 2015).To locate the optimal or near optimal solutions a heuristic algorithm such as Tabu search was developed (He et al., 2005). Several studies investigated on identifying the effectiveness of turn off /on on the machine tool energy consumption (Salido et al., 2013 ;Dai et al., 2013); Several authors have proposed mathematical programming and meta-heuristic way of solving bi-objective and multi-objective functions includes minimization of the total energy consumption and tardiness (Che et al., 2017; Zhang & Chiong, 2016).The existing research also demonstrates that minimization of maintenance energy costs, minimization of production energy costs and minimization of maximum completion time was modelled independently and solved(Song et al., 2014). Conversely, bi-objective functions such as completion time and energy consumption was modelled and solved using Non-dominated Sorting Genetic Algorithm II

*Contact: mallanjulap1@gmail.com*

(NSGA-II) in Flexible job shop problems (Yang et al., 2016a). Further three carbon efficiency indicators are put forward to estimate the carbon emission of parts and machine tools in the previous works and solved it through the NSGA- II coupled with local search algorithm based on neighbourhood search (C. Zhang et al., 2015).

Unequivocally, past literatures investigated major important objectives such as make-span, robustness paired up with the energy efficiency identify the relationship. Further concluded that there exists a positive correlation while maintaining balanced trade-off between energy consumption and robustness(Dai et al., 2013;Salido et al., 2013;Che et al., 2017;Yang et al., 2016b;Wu et al., 2013;Yang et al., 2016a)

Under the conditions of time-of-use tariffs a time-indexed integer programming formulation was developed that minimizes the carbon footprint and the cost of electricity. (H. Zhang et al., 2014).In-order to gain the targeted throughput along with the minimization of energy consumption, a control strategy was developed for a closed loop flowshop plant in the previous works(Mashaei & Lennartson, 2012). The investigation on energy utilization reduction was also explored by machine turning on and off. Alongside, energy consumption metric paired with cost of production was calculated for producing single product using analytical model (Zanoni et al., 2014). Further, a multi-objective mixed-integer programming (MIP) formulation including energy and completion time consideration on a single machine was done without considering the set-up times was bought into the picture (Fang et al., 2013). And later, a permutation flowshop problem with energy utilization constraints using MIP formulation was proposed, by considering both the discrete and continuous processing speed(H. Zhang et al., 2014).

In recent times, a significant progress has been made in bringing the CPU times lower in addressing an energy-efficient permutation flowshop scheduling problem while compromising the solutions(Yüksel et al., 2020b). On other hand an ensemble of meta-heuristics is proposed for the energy efficient blocking flowshop scheduling problem by taking small scale datasets into consideration (Kizilay et al., 2019). The imperative and necessitate for exploration in short term production planning through integrating energy in operations management has been demonstrated by range of authors (Gunasekaran & Ngai, 2012; Kleindorfer et al., 2005;Salgado & Pedrero, 2008;Gharbi et al., 2013a). Previous studies presented an energy scheduling problem which was formulated as a generalization of the cumulative scheduling problem, this itself an extension of the well-known parallel machine scheduling problem. The main characteristic of the problem solved is to meet the minimization of energy in shop floor (Artigues et al., 2013).

A mixed integer linear programming is formulated, defined by a set of periods where each one is characterized by a length, an electricity price, a maximal allowed power and an external demand by minimizing the total energy consumption and makespan in flexible manufacturing systems (FMSs) (Masmoudi et al., 2016; X. Li et al., 2017). Previous research also demonstrated an approach to minimize the make span and energy consumption of the dynamic scheduling problem for a flexible flowshop scheduling exhaustively. A dynamic flexible flowshop scheduling which is considered to be NP-hard problem is solved by adapting the novel algorithm based on an improved swarm optimization for pareto optimal solution (Mansouri et al., 2016;Tang et al., 2016).

*Contact: mallanjulap1@gmail.com*

Recently a machine learning based memetic algorithms was also proposed to solve the permutation problem flowline scheduling problem(Wang & Tang, 2017). Objectives to minimize both makespan and carbon footprint were considered simultaneously, which was solved by a multiobjective teaching learning-based optimization algorithm. Furthermore, three carbon-foot print reduction strategies were employed to optimize the scheduling results (W. Lin et al., 2015). A lower bound was developed for the permutation flowshop scheduling with sequence independent setups times based on waiting time-based relaxation (Gharbi et al., 2013b). A branch and bound algorithm were developed based on the NEH heuristic to solve the permutation flowshop scheduling problem by considering the energy as an important objective function (W. Liu et al., 2017).

Meta-heuristic on the other hand can establish an optimal solution or near optimal solutions with acceptable time consumptions. In fact, they are widely used in the combinatorial optimization. Of them multi-objective meta-heuristic algorithm is the non-dominated sorting genetic algorithm II (NSGA-II) (Deb et al., 2002). As it was suggested that without proper formation of a problem all algorithms will perform no better than random blind search. An ant colony optimization approach was also proposed in the previous works to solve the permutation flowshop scheduling problem with the possibility of outsourcing allowance to certain jobs (Neto & Godinho Filho, 2011). The three prominent metaheuristics algorithms namely simulated annealing, genetic algorithm, tabu search were also presented to solve the non-permutation flowline manufacturing cell with sequence dependent family setup times, thereby decreasing the computational time(S.-W. Lin et al., 2009). An enhanced co-evolutionary algorithm was proposed to resolve the multi-objective energy efficient task scheduling problem on a green data centre partially powered by the renewable energy(Lei et al., 2016). A summary of the selected literature survey is presented Table 1 along with objective function and algorithm used.

**Table 1:** Literature summary

*Contact: mallanjulap1@gmail.com*

| Objective Function | Algorithm used to Solve | Reference |
|---|---|---|
| Energy consumption, Total tardiness | Analytical Hierarchy Process | (Mouzon & Yildirim, 2008) |
| Energy efficiency | Co-evolutionary algorithm | (S.-W. Lin et al., 2009) |
| Energy consumption | Mathematical programming | (Fang et al., 2013) |
| Energy consumption | Branch and bound | (G.-S. Liu et al., 2013) |
| Energy efficiency | Genetic-simulated annealing | (Dai et al., 2013) |
| Energy consumption | Heuristic algorithm | (Mashaei et al., 2013) |
| Energy cost | Constraint Programming | (Mashaei et al., 2013) |
| Energy cost, $CO_2$ reduction | Mathematical programming-CPLEX | (H. Zhang et al., 2014) |
| Total Energy consumption cost | Genetic Algorithm | (Shrouf et al., 2014) |
| Makespan with restriction on Peak power | Heuristics | (Zanoni et al., 2014) |
| Energy consumption | Mathematical Programming – CPLEX | (Che et al., 2015) |
| Makespan, carbon footprint | Teaching learning-based algorithm | (W. Lin et al., 2015) |
| Energy consumption throughput | Heuristics | (C. Zhang et al., 2015) |
| Energy Efficiency and Productivity | Constraint Programming | (Yang et al., 2016b) |
| Energy efficiency | Branch and bound, Heuristics | (Lei et al., 2016) |
| Makespan, Energy | Heuristics | (Afshin Mansouri & Aktas, 2016) |
| Makespan, Energy consumption | Multi-objective genetic algorithm | (Mansouri et al., 2016) |
| Energy efficiency | Particle swarm optimization | (Tang et al., 2016) |
| Carbon efficiency | NEH-insertion algorithm | (Ding et al., 2016) |
| Energy consumption, tardiness | Local search | (Che et al., 2017) |
| Makespan, energy consumption | Efficient multi-objective heuristic algorithm | (Li et al., 2018) |
| Energy efficiency | Genetic algorithm | (Kizilay et al., 2019) |
| Makespan, energy consumption | Distributed algorithms | (Amiri & Behnamian, 2020) |
| Total tardiness, energy consumption | Heuristic | (Yüksel et al., 2020b) |
| Makespan, carbon emission | Multi-objective Iterated greedy | (Schulz et al., 2022) |

## 3. Energy Saving Scheduling Model for the PFLSP

As far as computer numerical control (CNC) machine is concerned, it has multiple energy source system (or subsystem).It can be classified into systems which consumes energy based on load known to be systems related to load energy consumption (SRLEC) and systems which does not consumes energy based on the load known as systems unrelated to load energy consumption (SULEC).The systems focussed on this dichotomy of SRLEC and SULEC is mainly based on the following characteristics as shown in the Figure 1. Sets, indices, parameters and variables used in the mathematical formulation is presented in equation-1.

*Contact: mallanjulap1@gmail.com*

**3.1 Electrical Components in Machine tool**

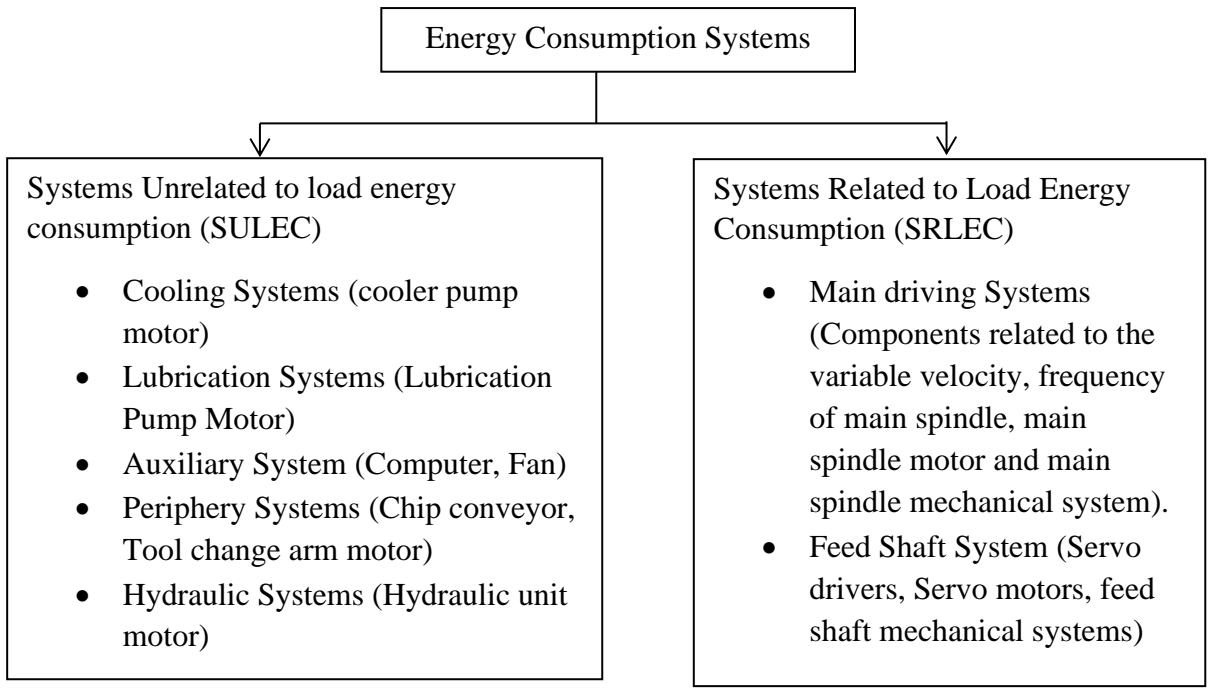

Figure 1: Energy consumption systems classification in real time

*Contact: mallanjulap1@gmail.com*

**Table 2:** Sets, indices, variables and parameters of mathematical model.

| | |
|---|---|
| Indices | |
| $i$ | Index for jobs; i=1,2…. n |
| $j$ | Index for machines j=1,2…m |
| $sp$ | Index for spindle speed |
| $f$ | Index for feed shaft |
| Parameters | |
| $n$ | number of jobs |
| $m$ | number of machines |
| $k$ | number of systems un-related to load energy consumptions |
| $P_{sp}$ | Power of main spindle speed |
| $f_{sp}$ | Feed speed |
| $M_{sp}$ | Friction |
| $B_{sp}$ | Coefficient of damping |
| $V_{sp}$ | Speed of cutting |
| $A_{sp}$ | Depth of cutting |
| $P_{axf}$ | Power of the feed shaft |
| $\Omega$ | Speed of the motor |
| $M_f$ | Motor torque due to friction |
| $B_f$ | Damping co-efficient of feed shaft |
| $T_f$ | Cutting torque |
| Positive variables | |
| $\pi_i$ | $i = 1,2,..,n$; is schedule |
| $t(i,j)$ | is the processing time when job i is processed on machine j |
| $C(\pi_i, j)$ | is the completion time when job $\pi_i$ is processed machine j |
| $FT$ | The total flowtime of jobs on machines |
| $EC$ | Total energy consumption |
| Binary variables | |
| $g_i(t)$ | 1 the i[th] subsystem is running |
| | 0 the i[th] subsystem is not running |

$$EC(t) = \sum_{i=1}^{n} \int_{t_0}^{t} g_i(t).C_i dt + P_{sp}(t, n, M_{sp}, B_{sp}, V_{sp}, A_{sp}, f_{sp}) dt$$
$$+ \sum_{f=1}^{m} \int_{t_0}^{t} P_{ax_f}(t, \Omega, M_f, B_f, f, T_f,) dt \qquad (1)$$

*Contact: mallanjulap1@gmail.com*

The equation (1) describes the energy consumption in the systems unrelated to the load energy consumption. Where $\sum_{i=1}^{n} \int_{t_0}^{t} g_i(t). C_i dt$ represents the energy of SULEC, n represents the number of SULEC. $g_i(t)$ denotes the used state of SULEC. The i$^{th}$ SULEC system power is represented as $C_i$ and it is independent of time t. The power of the main spindle system is represented as $P_{sp}(t, n, M_{sp}, B_{sp}, V_{sp}, A_{sp}, f_{sp})$, which is the function of friction $M_{sp}$, depth of cutting $A_{sp}$, rotating speed n, damping co-efficient $B_{sp}$, speed of cutting $V_{sp}$, and feed speed $f_{sp}$. The power of the feed shaft is given as $P_{ax_f}(t, \Omega, M_f, B_f, f, T_f,)$ which is the function of time t, speed of the motor $\Omega$, torque due to friction $M_f$, coefficient of damping $B_f$, feed speed $F$, and cutting torque $T_f$.

The SULEC state is represented by the equation (1) which reveals the state of energy consumption of the machine tool, when the machine is in the state of functioning. And it is given as $E_{in}(t) = \sum_{i=1}^{n} \int_{t_0}^{t} g_i(t). C_i dt$. The activated machine components of the machine which are ready for the operations is given by the fixed power $P_F = \sum_{i=1}^{n} g_i(t). C_i$. From the previous research (G.-S. Liu et al., 2013); (S. Li et al., 2018) it is found that, the author has pointed out the five states of machine tool includes cutting, air-cutting, starting, standby, idling. They are:

a) Cutting piece is represented by cutting state.
b) The retracting cutting tool or feeding is represented by the state of Air-cutting.
c) Idling state is represented by a certain speed of spindle's steady rotation state.
d) The speed and accelerating of the spindle motor is being represented through starting state.
e) Ensuring the operational readiness in the machine components and activating the machine tool is represented by the standby state

## 4. Problem definition

This study is motivated due to large number of organisations of late considering energy costs in their operations. In addition to this, researchers are also focusing on addressing flowshop scheduling with minimal energy consumption along with traditional measures like flowtime, makespan etc. The detailed energy saving model for the permutation flowshop scheduling problem is provided in Appendix. This study considers a general flowshop scheduling problem with "m" machines and "n" jobs. Wherein all the jobs are available at the same time and is defined as follows:

Set of jobs is denoted as J= *{1, 2, 3, ………, n jobs}* and set of machines as M=*{1, 2, . . . ,m machines}*. Every job passes through each machine only once and every job follows the same order. The assumptions of the considered permutation flowshop scheduling problem (PFLSP) are the following:

1. At most one operation can process on one machine at a time.
2. For the first machine, all jobs are available at time t=0.
3. In-between the jobs, there are no relationships precedented.
4. For every job between operations there exists precedence relationships.

*Contact: mallanjulap1@gmail.com*

5. No interruption is allowed once the operation starts on a machine until it is completed i.e., pre-emption is not permitted.
6. On every machine the sequences of jobs are consistent. i.e., to say, on machine 1 if the job is at the i$^{th}$ operation, then this can be concluded as on all the machines this job is at the i$^{th}$ position

The target of PFLSP is to identify a sequence of jobs that minimizes the total energy consumption and the total flowtime .

### 4.1. Model for Energy-saving based on the energy consumption.

Minimization of EC in flowshop scheduling can be achieved by minimizing the idle time of machine i.e., in-other words maximizing the occupancy of the machines by suitable schedule. Thereby decreasing the unrelated load energy consumption of energy systems, which means decreasing the energy consumption (EC). This section demonstrates the model for energy-saving based on the multi-objective scenario for minimizing the flowtime and energy consumption. The suggested model minimizes the total flowtime (FT) and energy consumption (EC).

The objective functions of the proposed model as follows:

$$Y = min\{F_{total}, E_F\} \tag{1}$$

The computed total energy consumption and total flowtime is given by the following formulae.

Completion time $C(\pi_i, j)$
$$C(\pi_1, 1) = t(\pi_1, 1)$$
$$C(\pi_i, 1) = C(\pi_{i-1}, 1) + t(\pi_i, 1) \quad i = 2, \ldots, n$$
$$C(\pi_1, j) = C(\pi_1, j-1) + t(\pi_1, j) \quad j = 2, \ldots, m$$
$$C(\pi_i, j) = max\{C(\pi_{i-1}, j), C(\pi_i, j-1)\}$$
$$+ t(\pi_i, j) \quad i = 2, \ldots, n; \ j = 2, \ldots, m$$

Total flowtime $F_{total}$
$$F_{total} = \sum_{i=1}^{n} C(\pi_i, m) \tag{1.1}$$

Stand by time of machine tools $T(\pi_i, j)$
$$T(\pi_i, 1) = 0 \quad i = 1, \ldots, n$$
$$T(\pi_1, j) = C(\pi_1, j-1) \quad j = 2, \ldots, m$$
$$T(\pi_i, j) = max\{C(\pi_i, j-1) - C(\pi_{i-1}, j), 0\}$$
$$i = 2, \ldots, n; j = 2, \ldots, m$$

Total Energy Consumption $E_F$
$$E_F = \sum_{j=1}^{m}\sum_{i=1}^{n} P_j^F . T(\pi_i\ j) \quad i = 1, \ldots, n; j = 1, \ldots, m \tag{1.2}$$

*Contact: mallanjulap1@gmail.com*

The equation (2.1) represents the total flowtime of $\pi_i$ $(where\ i \in 1,2,..,n)$ schedule when processed on $m$ machines. Equation (2.2) represents the total energy consumption of $\pi_i$ $(where\ i \in 1,2,..,n)$ schedule when processed on m machines which is essentially the unutilized energy due to either job or machine being idle when the previous job or current machine under the operation.

## 5. Multi-objective evolutionary algorithm for Permutation flowshop scheduling Problem

Based on the equation 2.1 and 2.2 this study is subjected to two main objective functions i.e., flowtime and energy consumption. Since the problem considered can be classified as NP-hard problem due to its computational complexity in the non-deterministic polynomial time with its bi-objectives(G.-S. Liu et al., 2013). Albeit there exist various approaches to address this problem, this study employed metaheuristic approach. The aim of these algorithms is to explore the search space in order to find the approximate pareto solutions. This study used NSGA-II which is a meta-heuristic algorithm for solving these multi-objective problems.

### 5.1. Proposition of NSGA-II algorithm

The NSGA-II algorithm ameliorate the selection operation on the grounds of genetic algorithm in-order to fathom the multi-objective optimization problem. The NSGA-II algorithm can provide the Pareto frontier solution by incorporating the elite retention strategy and the congestion comparison strategy into the selection operation(Deb et al., 2002)

a. Non-dominated sorting method

In the proposed optimization problem with $O$ fitness functions, there exists a $\phi_h(x_A) \leq \phi_h(x_B)$ $\phi_h$ is the objective value on the $h$ optimization objective under the hypothesis of $\forall h \in \{1,2,\ldots,O\}$. This can be demonstrated as individual "A" dominating individual "B". An individual solution is known to be non-dominated one if there exists no individual solution in a given population of "A". Further, based on the degree of inferiority or superiority each individual received the rank label.

b. Comparison operator for crowding distance criteria

In the sequence of selecting the genes based on the equal non-dominated level crowding distance criteria is being used. The crowding distance of the genes at the lower bound and upper bound is given as infinity. This can be demonstrated as follows:

$$d_m = \sum_{h=1}^{O} |f_h^{m+1} - f_h^{m-1}| \qquad (2)$$

From the above formula (2), $f_h^{m+1}$ and $f_h^{m-1}$ are considered to be the adjoining individual fitness function values of individual $m$ and $d_m$ is the crowding distance of individual.

c. Strategy for elite retention

As demonstrated in Figure 2. The fundamental aspects of the retention strategy are to retain the best individuals in the offspring from the previous(parent) generation. For the given population

*Contact: mallanjulap1@gmail.com*

size of N, unifying the parent population $P_t$s with the subpopulation $Q_t$ to generate the overall population as $R_t$, where t = 1..., N. Hence, the resulting population size is 2N. These 2N individuals are further ranked by non-dominated sorting.

5.2. Steps in Optimization model

Based on the NSGA-II for multi-objective optimization problem, the procedure for solving the proposed model, which is shown in Figure 2, is summarized as follows:

Step 1: Set initial algorithm parameters such as crossover and mutation probabilities at counter t=0, further generate the random schedules in which each gene represents a schedule P(t).

Step 2: Based on the dual important criterions namely ranking and crowding distance, every individual from the population set P(t) are arranged according to the non-dominated sorting algorithm. The arrangement of this original population set P(t) is done through the proposed energy saving model.

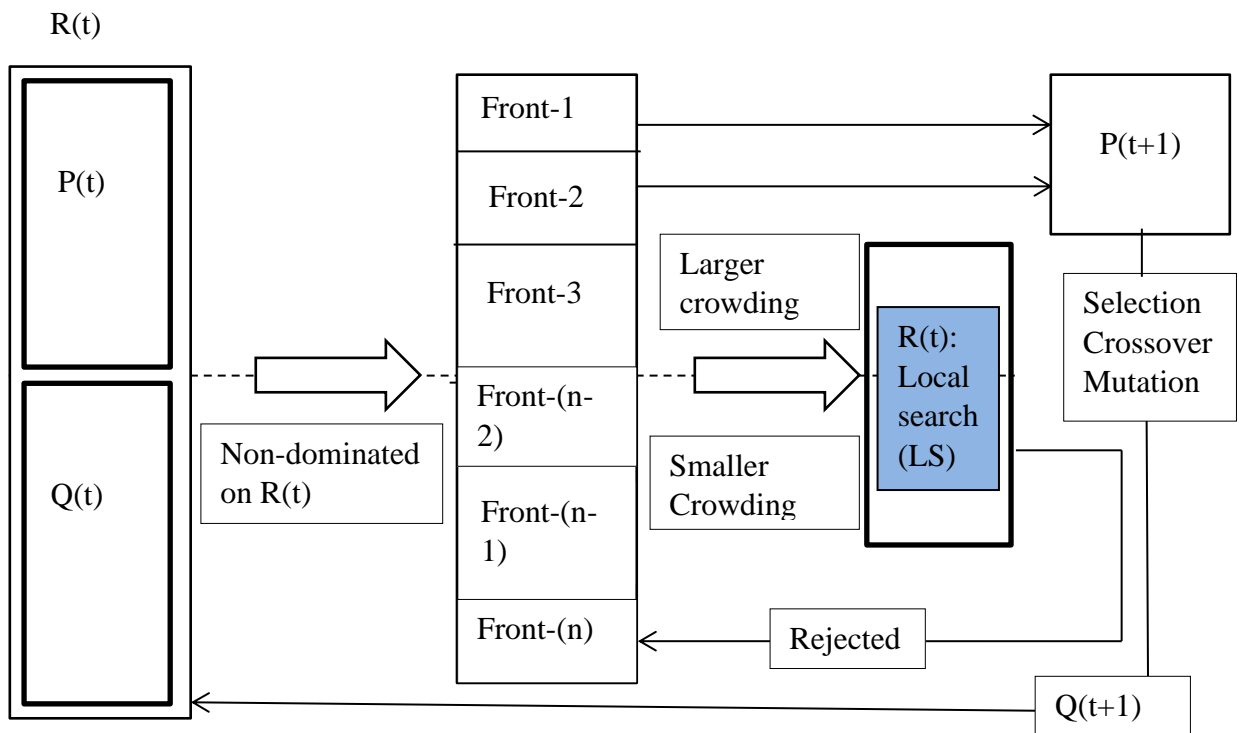

Figure 2: Diagrammatical view of process flowchart for the proposed NSGA-II algorithm.

Step 3: The parent population Q(t), is formed based on the tournament selection on the P(t) where genes with larger crowding distance value have significant probability to be selected.

Step 4: Orchestrate the internal mechanism of operators (crossover and mutation) on parent population Q(t), in order to generate the offspring population Q(t). Consequently, Q(t) population undergo fitness evaluation for each gene.

Step 5: Unifying population to form R(t) = P(t)UQ(t), and carry out non-dominated sorting on R(t).

*Contact: mallanjulap1@gmail.com*

Step 6: Post-unification the solutions from non-dominated sorting with Rank 1, undergo local search to identify the better solutions that may exist in the neighbourhood of original non-dominated sorted solutions before next iteration. In-order to define the neighbourhood of the non-dominated solution local search operators are designed. To replace the current non-dominated sorted solutions, variable neighbourhood descent algorithm is implemented to obtain a better non-dominated solution. With the proposed local search, best solutions around the neighbourhood of non-dominated sorted solutions can be selected to fast-track convergence

Step 7: Based on the maximum number of pre-defined iterations, the algorithm selects whether to proceed to next iteration i.e., t=t+1, by jumping back to Step-3 or to terminate iterative search with pareto fronts of the planning solutions as outcome.

### 5.1. Local search operators and Pseudo code:

The intermediate local search (LS) technique that is implanted with in the NSGA-II consists of three search operators in addition to classical crossover and mutation operators. These three LS operators are developed to find the nearest neighbourhood solutions in the non-dominated sorted solutions. Initially, the first LS operates three different types of operations essentially includes the 1. swapping the sequence of jobs of the initial solution, 2. Performing the reversion operator and finally 3. Performing the neighbourhood operation. Thus P(t) to get two additional solutions in each operator except for the operator 3. Line-6 in the pseudo code represents the same. Further, post these local search operations, the solutions in set P undergo non-dominated sorting to get Rank 1 solution set $N_1$. This further ensures that no solution from the set P can dominate any solution in set $N_1$. These sequences of operations are represented from lines 8-21. With the NSGA-II crowding distance criteria, the solutions from $N_1$ set undergo crowding distance calculation. The solution with the highest crowding distance (X) is selected to be the new solution. Finally, the identified new solution $X$ is compared with the best solution $X_{best}$. If $X$ is same as $X_{best}$ then algorithm chooses the next immediate LS operator. This is represented from line 29 and 30. Each LS operator generates two more neighbour solutions except operator 3 which typically generates 10 solutions. If the $X$ is not the same $X_{best}$ then the best fitness chromosome in the neighbourhood is better than ($X_{best} \Leftarrow X$) $X$ is assigned to be the best. This is represented in the line 26. Terminate the loop if all the LS operators cannot improve the current solution.

1. Swapping of sequence operator.

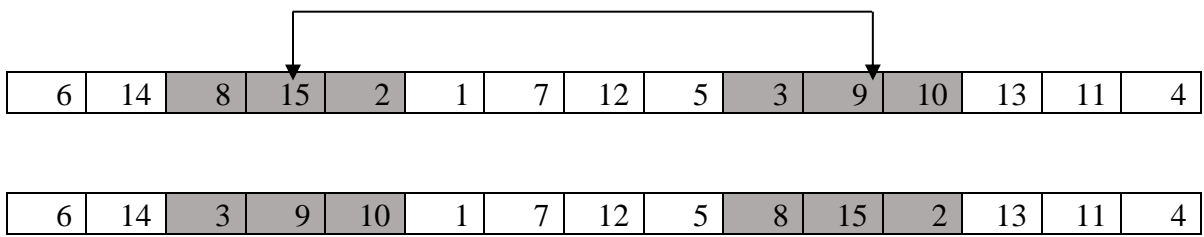

*Contact: mallanjulap1@gmail.com*

2. Reversion operator.

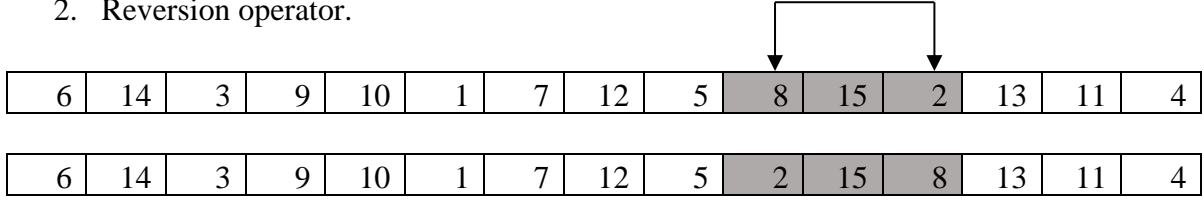

3. Neighborhood operator.

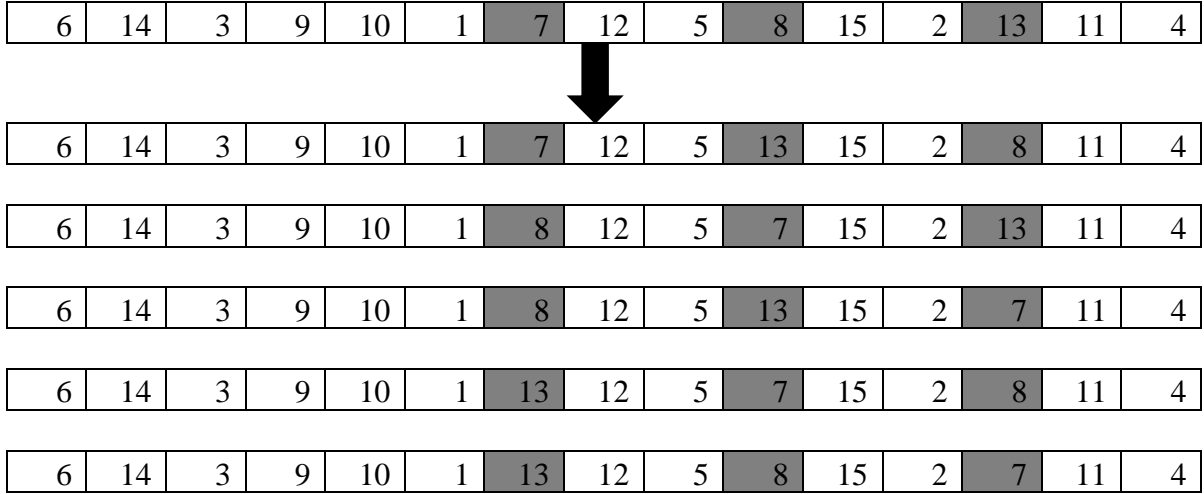

Figure 3: Local search operators of NSGA-II-LS algorithm.

**Input**: $A\ set\ of\ search\ operators\ O_a(a = 0, 1, 2) for\ finding\ neighbor\ solutions;$
$The\ three\ objective\ functions\ EC(x), FT(x); Maximum\ generation\ G_{max};$
$An\ initial\ solution\ X_0; An\ initial\ action\ a = 0;$
**Output**: $An\ optimal\ solution\ X_{best};$
1: $X_{best} \Leftarrow X_0$
2: $a_p \Leftarrow 0\ //\ a_p\ is\ the\ previous\ action$
3: flag $\Leftarrow 0$
4: $g \Leftarrow 1$
5: While g < $G_{max}$; do
6: P $\Leftarrow O_a(X_0) + X_0$
7: $N_1 \Leftarrow \emptyset$
8: for each solution p∈P do
9: $S_p \Leftarrow \varphi\ //\ S_p\ is\ used\ to\ store\ the\ solutions\ that\ are\ dominated\ by\ solution\ p$
10: $n_p \Leftarrow 0\ //\ n_p\ is\ the\ number\ of\ times\ that\ p\ is\ dominated$
11: $r_p \Leftarrow 0\ //\ r_p\ represents\ the\ rank\ of\ p$
12: for each solution q ∈ P do
13: $if\ p\ dominates\ q\ \cap p \neq q\ then\ S_p \Leftarrow S_p \cup \{q\}$
14: $end\ if$

15: $if\ q\ dominates\ p\ \cap p\ \neq q\ then\ n_p \Leftarrow n_p + 1$
16: $end\ if$
17: $end\ for$
18: $if\ n_p == 0\ then\ r_p \Leftarrow 1$
19: $N_1 \Leftarrow N_1 \cup \{p\}$
20: end if
21: end for
22: $if\ |N_1| == 1\ then\ X \Leftarrow N_1[0] //\ Set\ N_1\ has\ one\ element, let\ X\ be\ the\ element; else$
23: $for\ each\ solution\ p\ \in N_1\ //\ \text{calculation of the crowding distance for Cp of p and X} \Leftarrow \max(N_1)$
24: $end\ for$
25: $end\ if$
26: $if\ X\ dominates\ X_{best} \cap X \neq X_{best}\ then\ X_{best} \Leftarrow X$
27: $else\ a_p \Leftarrow a; \text{flag} + +$
28: $a \Leftarrow flag\%3; end - if$
29: $if\ a_p - a == 2\ then\ \text{break}; end - if$
30: $g + +; end - while$

## 6. Implementation of proposed NSGA-II_LS algorithm on (S. Li et al., 2018) **dataset**.

The proposed algorithm was coded in MATLAB programming language and run on an Intel Xenon CPU 3.50 GHz with 16.0 GB RAM under Windows 10 Enterprise. In order to validate the viability and efficacy of the planned method, the modified NSGA-II algorithm tested with the same numerical experiments reported in the paper for the multi-objective problems. The results obtained through the improved algorithm is compared with the existing literature solutions (S. Li et al., 2018).

To evaluate the proposed method, under different problem sizes, we need to consider flowshop scheduling problem for the preliminary experiment. In this numerical illustration, we considered the same example of input data (S. Li et al., 2018). There are 15 jobs and 5 machines. These 5 machines have energy consumption and processing time of each job is shown in Table 3.

*Contact: mallanjulap1@gmail.com*

**Table 3:** 15 Jobs x 5 machines fixed power and job processing time (min)

| Machines | M1 | M2 | M3 | M4 | M5 |
|---|---|---|---|---|---|
| Fixed Power(Whr) | 769 | 802 | 1290 | 967 | 1166 |
| Job number | | | | | |
| 1 | 3 | 4 | 6 | 10 | 3 |
| 2 | 4 | 5 | 2 | 8 | 8 |
| 3 | 7 | 10 | 8 | 4 | 7 |
| 4 | 9 | 10 | 2 | 2 | 6 |
| 5 | 2 | 2 | 5 | 9 | 9 |
| 6 | 2 | 1 | 1 | 8 | 3 |
| 7 | 5 | 7 | 8 | 2 | 5 |
| 8 | 2 | 9 | 2 | 9 | 8 |
| 9 | 9 | 7 | 3 | 8 | 1 |
| 10 | 8 | 5 | 7 | 2 | 2 |
| 11 | 9 | 6 | 9 | 4 | 7 |
| 12 | 7 | 9 | 3 | 2 | 4 |
| 13 | 8 | 8 | 2 | 2 | 9 |
| 14 | 1 | 2 | 6 | 5 | 9 |
| 15 | 8 | 2 | 10 | 1 | 4 |

### 6.1. Experimental setup through Taguchi methodology

It should also be noted that the performance of the compared algorithms is affected by algorithm's parameters settings. The preliminary experiments through design of experiments are conducted under a set of parameters to find the best combination. Parameters of the NSGA-II can be tuned through various methodologies through exhaustive experimentations. In this proposed study we adopted Taguchi methodology to run the selected experiments in-order to select best generation size, population size, crossover probability and mutation probability. Table-4 represents the various NSGA-II parameters along with the levels chosen to conduct the Taguchi analysis.

**Table 4:** NSGA-II parameters with various levels.

| | Level-1 | Level-2 | Level-3 | Level-4 |
|---|---|---|---|---|
| **Pop** | 25 | 50 | 100 | 200 |
| **Gen** | 10 | 25 | 50 | 100 |
| **Mutation** | 0.05 | 0.06 | 0.07 | 0.08 |
| **Crossover** | 0.5 | 0.6 | 0.7 | 0.8 |

With the adaption of Taguchi methodology, we can effectively reduce the experimental costs and obtain high quality results by fewer number of experiments through orthogonal arrays. The two main groups of Taguchi method i.e., control factors or design parameters and noise factor (hard to control) are designed to obtain the desired outcome. The deflection between the aimed values and experimental values in the Taguchi method is calculated through loss function. Usually, $L_a(b, c)$ is the representation of the orthogonal array. Here '$c$' indicates the number of factors and '$b$'

*Contact: mallanjulap1@gmail.com*

indicates the number of levels. In this study $L_{16}$ design based on the number of factors i.e., 4 and numbers of levels i.e., 4.

6.1.1 Experimental runs with Flowtime as response variable:

Table-4 represented the orthogonal arrays that we experimented with the response variable as flowtime. In NSGA-II for the best parameter selection with 4 degree of freedom, the mean S/N ratios are calculated in terms of mean response as shown in tables 5. Table-5 represents the response table for means with population size being ranked -1 followed by mutation, generation size and crossover. S/N ratio for mean response is presented in figure-4 and figure-5, helpful in selection of optimum combination of NSGA-II for the PFLSP with flowtime as the response variable. From the figure-4 it indicates that population size of 200, generation size of 50, mutation probability as 0.05 and with the crossover probability of (either 0.5, or 0.6 and 0.8) can be the best to be selected. Before we conclude these parameters, we conducted additional study with the response variable as energy consumption as the response variable.

**Table 4:** Taguchi Orthogonal $L_{16}$ array with flowtime as response variable.

| Gen | Pop | Crossover | Mutation | Flowtime |
|---|---|---|---|---|
| 10 | 25 | 0.5 | 0.05 | 912 |
| 10 | 50 | 0.6 | 0.06 | 913 |
| 10 | 100 | 0.7 | 0.07 | 909 |
| 10 | 200 | 0.8 | 0.08 | 921 |
| 25 | 25 | 0.6 | 0.07 | 920 |
| 25 | 50 | 0.5 | 0.08 | 915 |
| 25 | 100 | 0.8 | 0.05 | 917 |
| 25 | 200 | 0.7 | 0.06 | 917 |
| 50 | 25 | 0.7 | 0.08 | 916 |
| 50 | 50 | 0.8 | 0.07 | 916 |
| 50 | 100 | 0.5 | 0.06 | 916 |
| 50 | 200 | 0.6 | 0.05 | 924 |
| 100 | 25 | 0.8 | 0.06 | 912 |
| 100 | 50 | 0.7 | 0.05 | 923 |
| 100 | 100 | 0.6 | 0.08 | 909 |
| 100 | 200 | 0.5 | 0.07 | 923 |

*Contact: mallanjulap1@gmail.com*

**Table 5:** Flowtime response table for means

| Level | Gen | Pop | Crossover | Mutation |
|---|---|---|---|---|
| 1 | 913.8 | 915.0 | 916.5 | 919.0 |
| 2 | 917.3 | 916.8 | 916.5 | 914.5 |
| 3 | 918.0 | 912.8 | 916.3 | 917.0 |
| 4 | 916.8 | 921.3 | 916.5 | 915.3 |
| Delta | 4.3 | 8.5 | 0.3 | 4.5 |
| Rank | 3 | 1 | 4 | 2 |

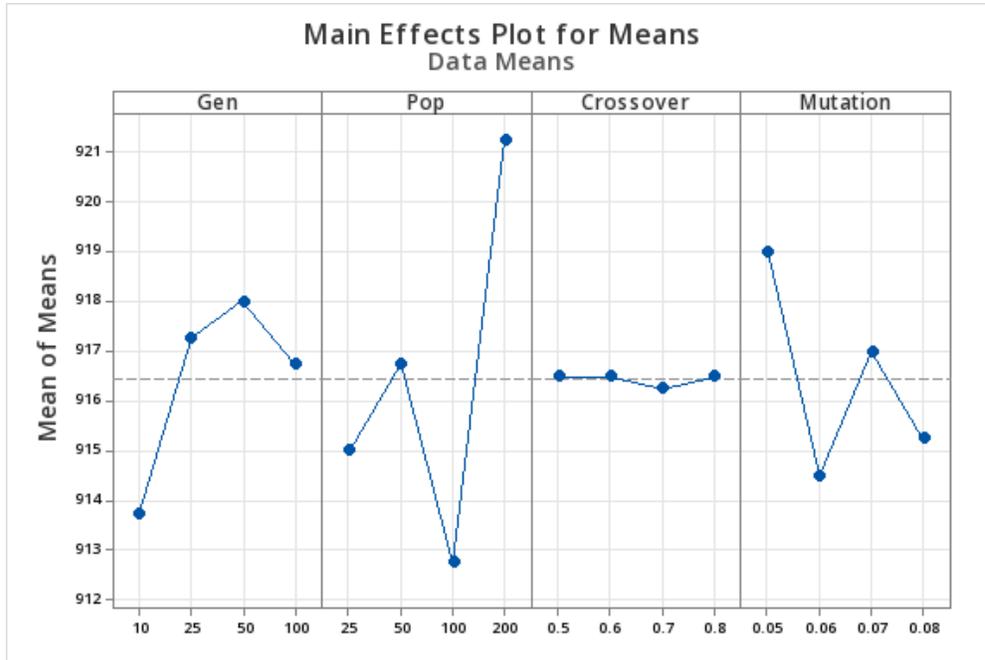

Figure 4: Effect of NSGA-II parameters on minimum flowtime response

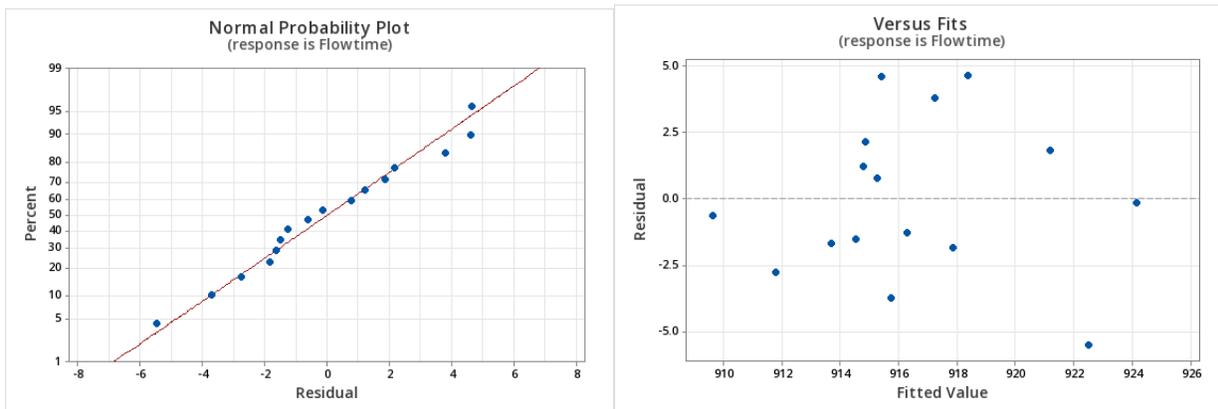

*Contact: mallanjulap1@gmail.com*

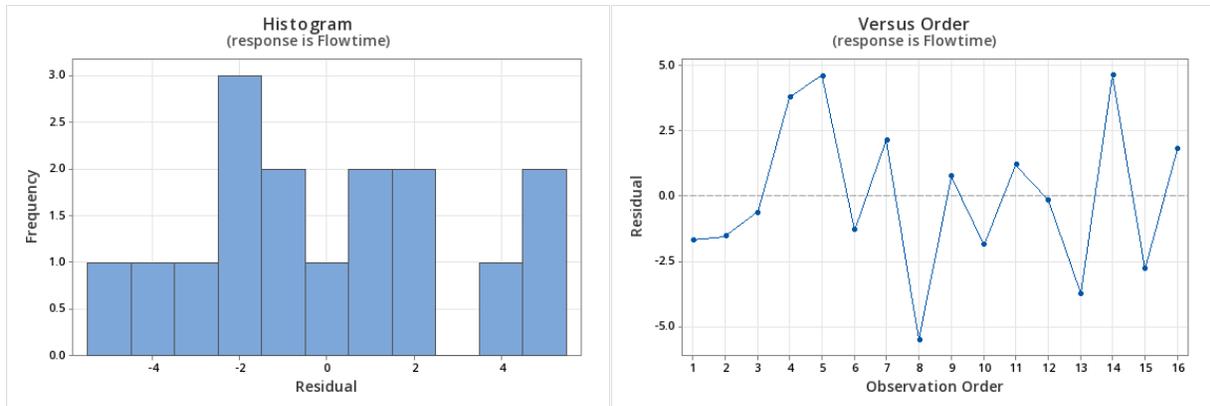

Figure 5: Residual plots for means to obtain minimum flowtime response.

6.1.2 Experimental runs with energy consumption as response variable:

From the table-4 we adopted the same number of experimental designs but in this table-6 analysis we altered the response variable to be the energy consumption (EC). Basically, we performed this additional analysis to finalize the parameter setups. Table-6 represents the orthogonal experimental design with EC response. From the table-7 represents the response table for means with mutation ranked -1 followed by crossover, generation size and population size.

**Table 6:** Taguchi Orthogonal $L_{16}$ array with EC as response variable.

| Gen | Pop | Crossover | Mutation | Energy Consumption |
|-----|-----|-----------|----------|--------------------|
| 10  | 25  | 0.5       | 0.05     | 1290.8             |
| 10  | 50  | 0.6       | 0.06     | 1207.8             |
| 10  | 100 | 0.7       | 0.07     | 1145.8             |
| 10  | 200 | 0.8       | 0.08     | 1207.8             |
| 25  | 25  | 0.6       | 0.07     | 1207.8             |
| 25  | 50  | 0.5       | 0.08     | 1207.5             |
| 25  | 100 | 0.8       | 0.05     | 1242.2             |
| 25  | 200 | 0.7       | 0.06     | 1290.4             |
| 50  | 25  | 0.7       | 0.08     | 1207.8             |
| 50  | 50  | 0.8       | 0.07     | 1209.8             |
| 50  | 100 | 0.5       | 0.06     | 1207.5             |
| 50  | 200 | 0.6       | 0.05     | 1348.7             |
| 100 | 25  | 0.8       | 0.06     | 1253.4             |
| 100 | 50  | 0.7       | 0.05     | 1334.2             |
| 100 | 100 | 0.6       | 0.08     | 1290.4             |
| 100 | 200 | 0.5       | 0.07     | 1145.8             |

**Table 7:** Energy consumption response table for means

| Level | Gen | Pop | Crossover | Mutation |
|---|---|---|---|---|
| 1 | 1213 | 1240 | 1213 | 1304 |
| 2 | 1237 | 1240 | 1264 | 1240 |
| 3 | 1243 | 1221 | 1245 | 1177 |
| 4 | 1256 | 1248 | 1228 | 1228 |
| Delta | 43 | 27 | 51 | 127 |
| Rank | 3 | 4 | 2 | 1 |

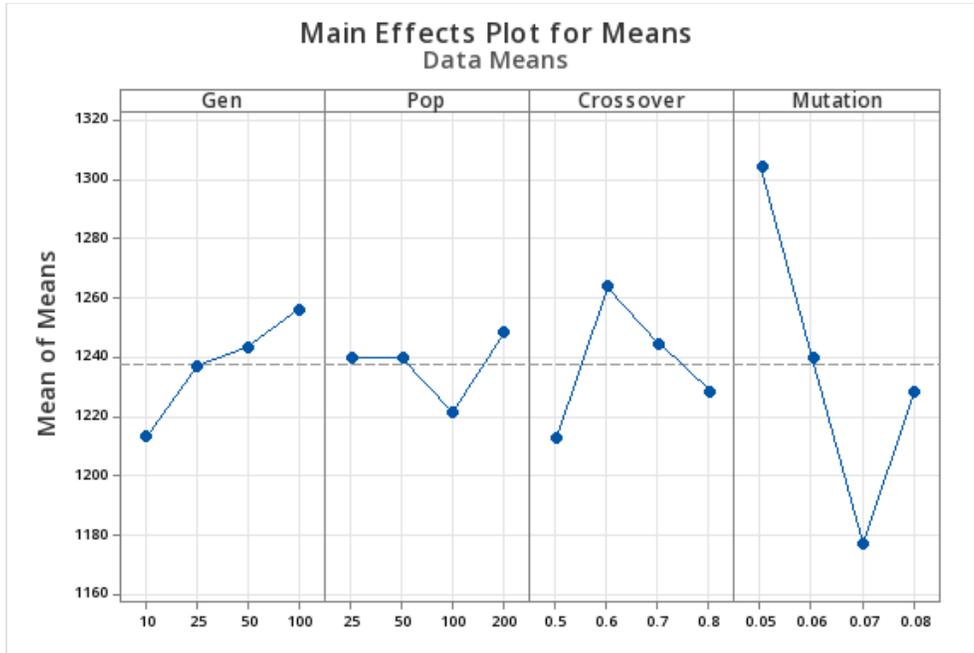

Figure 6: Effect of NSGA-II parameters on minimum energy consumption response

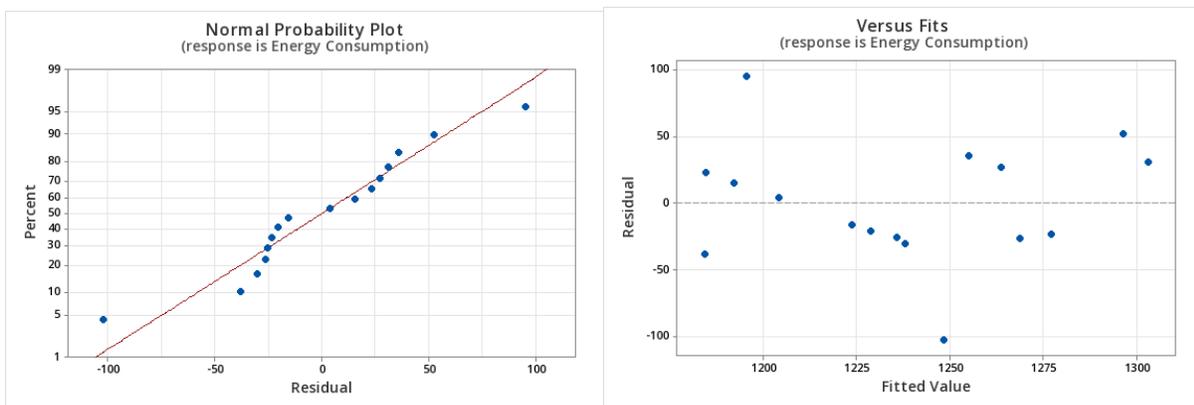

*Contact: mallanjulap1@gmail.com*

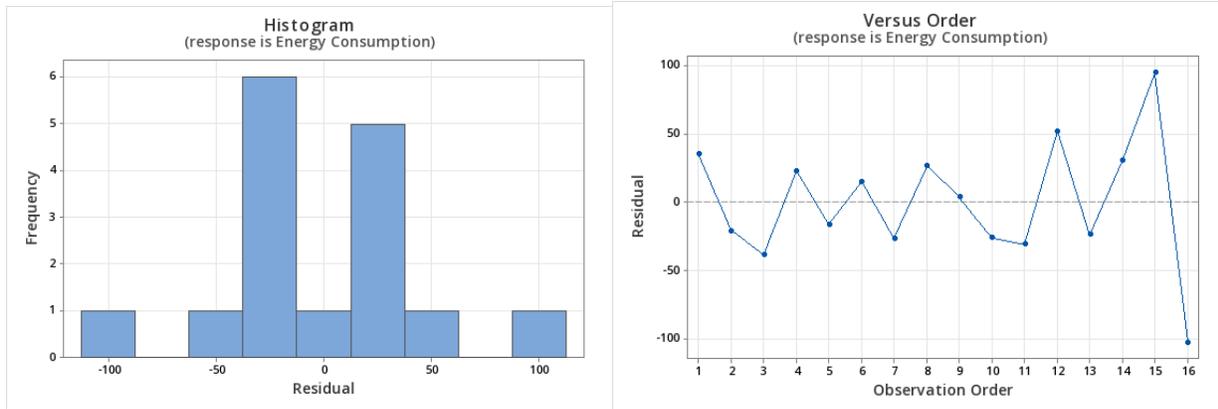

Figure 7: Residual plots for means to obtain minimum energy consumption response

From the figure-6 and figure-7 it indicates that population size of 200, generation size of 100, mutation probability of 0.05 and with the crossover probability of 0.6 we can achieve the minimum energy consumption. Based on the results from section 6.1.1 and 6.1.2 we selected the population size 200, generation size 50, mutation probability 0.05 and with the crossover probability 0.6 are selected to run the final run on the literature dataset i.e., from table-3. The results obtained through combinations of parameters are presented in table 4. The illustration of the preliminary experiment and the trade-off relationship between the total energy consumption and flowtime is demonstrated through the Pareto distribution chart in figure 3.

**Table 8:** The optimal solutions obtained at Crossover probability (Pc=0.6) and mutation probability (Pm=0.05)

| SI. No | Literature Solution by (S. Li et al., 2018) | | Proposed Solution | |
|---|---|---|---|---|
| | Total flowtime(min) | Energy consumption (Whr) | Total flowtime(min) | Energy consumption (Whr) |
| 1 | 912 | 1348.7 | 909 | 1348.7 |
| 2 | 915 | 1290.4 | 910 | 1309.8 |
| 3 | 919 | 1207.8 | 913 | 1290.4 |
| 4 | 919 | 1207.8 | 915 | 1290.4 |
| 5 | 919 | 1207.8 | 916 | 1207.8 |
| 6 | 932 | 1145.8 | 932 | 1145.8 |

*Contact: mallanjulap1@gmail.com*

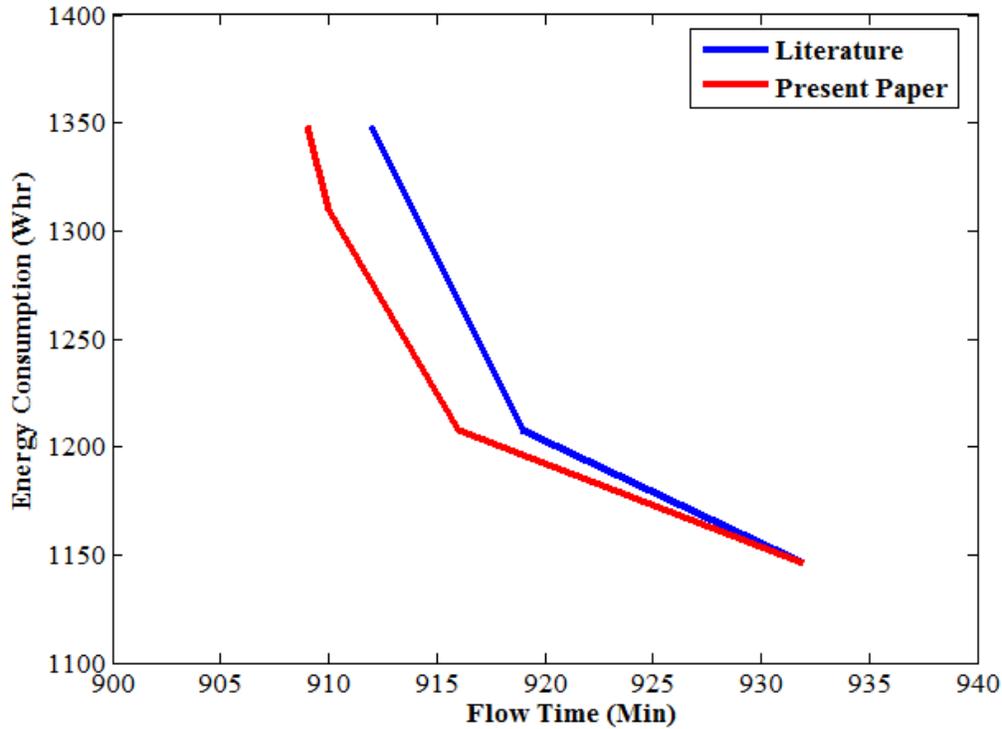

**Figure 8:** The comparative approximate pareto optimal front distribution of total FT versus EC

### 6.2. Implementation and performance evaluation of proposed algorithm.

The proposed NSGA-II_LS algorithm is implemented and evaluated using the results published (S. Li et al., 2018). To prove the effectiveness of algorithm we conducted numerical experiments with one (Taillard, 1993) benchmark dataset. It includes 9 different job sets (Ta20, Ta50, Ta100) processed on (5,10,20) machines as presented in Table 7. The time of processing of these jobs satisfy the uniform distribution of UD[1, 99], and the fixed power consumption satisfy the uniform distribution of UFP[700 1500]. The proposed improved algorithm was run for 10 times and the best combination of flowtime and its fixed power consumption is given. Energy consumption values considered in the problem datasets are presented in table 9.

**Table 9:** Energy consumption (Whr) for all the machines considered in Taillard datasets

| Machines | Fixed Power (Whr) |
|---|---|
| M1 | 769 |
| M2 | 802 |
| M3 | 1290 |
| M4 | 967 |
| M5 | 1166 |
| M6 | 1003 |
| M7 | 1211 |

*Contact: mallanjulap1@gmail.com*

| M8  | 1321 |
| M9  | 989  |
| M10 | 1411 |
| M11 | 782  |
| M12 | 980  |
| M13 | 1005 |
| M14 | 1333 |
| M15 | 867  |
| M16 | 1209 |
| M17 | 781  |
| M18 | 809  |
| M19 | 1113 |
| M20 | 977  |

The results obtained to the benchmark problems for the sample instance through the experiments were presented in the table 10 along with the literature solutions(S. Li et al., 2018), where FT1 signifies the flowtime minimum and EC1 signifies the matching energy consumption. EC2 signifies the energy consumption minimum and FT2 signifies the matching flowtime. According to the results obtained through the benchmark problems for one dataset, it can be interpreted that as the problem size increases the flowtime and energy consumption increases. Besides, different to the flowtime minimum scheduling, the schema for the most energy saving schedule can achieve by decreasing rate of energy consumption i.e., by computing $((E1-E2)/E1)$ with the total FT increasing by $((T2-T1)/T1)$. Out of all the experiments conducted on the benchmark problems, the proposed algorithm performed well in minimizing the flowtime and minimizing the energy consumption when compared to (S. Li et al., 2018).

### 6.3. Comparison of NSGA-II_LS with classical NSGA-II algorithm.

Fundamentally, we described the improvement as "Lower the percent FT increase along with the higher the percent EC increase is better". The proposed algorithm performed better on both the objectives i.e. FT and EC minimization in 5 out of 9 cases. On FT objective our algorithm performed better in 8 out of 9 cases. On EC objective our algorithm performed better in 5 out of 9 cases. Overall, we achieved 47.05% and 15.44% improvement in FT and EC minimization respectively. Whereas the literature solutions the average percentage difference in FT and EC were found to be 3.74% and 12.26% respectively. The proposed algorithm performed well with the overall improvement in the total flowtime and energy consumption. Thus, with the less total FT increase, the energy-efficient scheduling for the permutation flowline can decrease the energy consumption.

*Contact: mallanjulap1@gmail.com*

**Table 10:** The comparative results of proposed NSGA-II_LS algorithm and the literature solution.

|  | Best sol-FT |  | Best Sol-EC |  | % Diff (Prop. Solution)-NSGA-II_LS |  | % Diff (Lit Solution*)-NSGA-II |  |
|---|---|---|---|---|---|---|---|---|
| Problem | FT-Best-1 | EC-1 | FT-2 | EC-best-2 | (FT2-FT1)/FT1 | (EC1-EC2)/EC1 | (FT2-FT1)/FT1 | (EC1-EC2)/EC1 |
| Ta20x5 | 14502 | 13890 | 14650 | 12433 | **1.02** | 10.49 | 1.4 | **13.2** |
| Ta20x10 | 23757 | 86507 | 24212 | 72129 | **1.91** | 16.6 | 6.6 | **19.5** |
| Ta20x20 | 34988 | 293664 | 35839 | 278130 | **2.43** | 5.28 | 6.3 | **7.8** |
| Ta50x5 | 76690 | 26364 | 78417 | 20153 | **2.25** | **23.56** | 2.7 | 8 |
| Ta50x10 | 100650 | 96926 | 104030 | 71084 | 3.36 | **26.66** | **3** | 20.7 |
| Ta50x20 | 138280 | 346350 | 141150 | 313390 | **2.08** | **9.52** | 2.5 | 7.6 |
| Ta100x5 | 297390 | 24759 | 299880 | 23514 | **0.84** | 5.03 | 3.1 | **10.1** |
| Ta100x10 | 355213 | 141423 | 361950 | 111160 | **1.89** | 21.3 | 4.6 | 11.5 |
| Ta100x20 | 448923 | 573762 | 458221 | 504452 | **2.07** | 12.07 | **3.5** | 11.9 |
|  |  |  |  | Avg | 1.98 | 14.50 | 3.74 | 12.26 |
|  |  | % Improvement |  |  | 47.05 | 15.44 |  |  |

In this study, the NSGA-II_LS is compared with the classical NSGA-II in 10 different datasets 10 times, available from the literature (S. Li et al., 2018). In-order to assess the efficiency of the NSGA-II_LS algorithm, the true pareto front of the considered problems may be unknown. Therefore, for a given instance of problem, we treated every independent run non-dominated solution for these two different algorithms as pareto front for that instance. The best solutions are highlighted in bold font. A statistical test provides a significant comparison different due to the nature of the problem is stochastic Therefore, we conducted the Wilcoxon sign rank test to gauge the results generated by two different algorithms. We maintained 95% (corresponding to α = 0.05) confidence interval for all the tests. The main reasons for the NSGA-II_LS algorithm performance is superior to NSGA-II is because of the following: First, the parameter selection through Taguchi methodology diversified the population and further enhanced the solution quality due to the fact that NSGA-II is highly subtle to these parameters. Second, the local search technique with swapping of sequence operator, reversion operator and neighborhood operator identified the dual population nearby the originally selected population by classical NSGA-II. Overall, the neighborhood strategy has a positive effect on the selected population from non-dominated sorted solutions.

### 6.2. Insights from the implementation of proposed algorithm on Taillard datasets.

In this section we present the results obtained by the proposed NSGA-II_LS for the datasets selected from Taillard (1993). These datasets include i.e., (Ta20X5), (Ta20X10), (Ta20X20), (Ta50X5), (Ta50X10), (Ta50X20), (Ta 100X5), (Ta 100X10), (Ta100X20). In every dataset in-total we have investigated 10 sets. This is equivalent to (9x10) 90. We ran the algorithm for 10 times in solving every dataset i.e. (90x10) 900 instances. The best results obtained for each dataset in terms of FT and EC are presented in Table 11, Table 12 and Table 13. Since there are no standard solutions in the literature however, we present the results obtained by proposed algorithm in-terms of FT and EC objectives. This may help future studies for further evaluation. The following conclusions can be drawn.

*Contact: mallanjulap1@gmail.com*

1. As described in the section 6.1 the improvement is defined in the published literature (S. Li et al., 2018) as "Lower the percent FT increment along with the higher the percent EC increment is better". From the Table 11, Table 12 and Table 13 the results indicate that as the number of machines increase across Ta20, Ta50 and Ta100 jobs the average percentage increment in FT decreased from 3.75% to 3.09%. On other hand, the average percentage decrease of EC increased significantly from 11.41% to 16.89%.
2. Further, from the table-11, table-12, and table-13 the results indicate that as the number of jobs increase from Ta20 to Ta100 on any given number of machines (per say 5 machines) the average percentage of FT dropped from 5.14% to 4.87%. Although, this trend is not observed on 10 machines but this trend is observed when Ta20, Ta50, Ta100 jobs processed on 20 machines. On other hand, the average percentage decrease in EC increased from 14.90% to 27.14% when Ta50, Ta50 and Ta100 jobs processed on 5 machines. Similar trend is observed when these jobs processed on 20 machines except on 10 machines.
3. From the table-11 shows the results obtained for (Ta20X5), (Ta20X10), (Ta20X20). The average difference between FT and EC is calculated for each instance as well as each size of the dataset. It is observed that the difference in FT's is varying between 0.16% to 10.47% with an average of 3.47%. Similarly, EC's is varying between 2.87% to 20.99% with an average of 11.41%.
4. Table-12 shows the results obtained for (Ta50X5), (Ta50X10), (Ta50X20). The average difference between FT and EC is calculated for each instance as well as each size of the dataset. It is observed that the difference in FT's is varying between 0.07% to 6.83% with an average of 3.50%. Similarly, EC's is varying between 2.07% to 26.66% with an average of 14.71%.
5. Further, from the table-13 it is observed that for (Ta100X5), (Ta100X10), (Ta100X20). The average difference between FT and EC is calculated for each instance as well as each size of the dataset. It is observed that the difference in FT's is varying between 0.35% to 11.66% with an average of 3.09%. Similarly, EC's is varying between 2.84% to 43.73%% with an average of 16.89%.

*Contact: mallanjulap1@gmail.com*

**Table 11:** Taillard 20 jobs datasets with various machines configurations and their corresponding solutions with flowtime and energy consumption.

| Problem size (jobs x m/c) | Dataset | Total flowtime | | Energy consumption | | Percentage variation | |
|---|---|---|---|---|---|---|---|
| | S. No | Min (FT1) | EC1 | FT2 | Min (EC2) | (FT2-FT1)/FT1 | (EC1-EC2)/EC1 |
| Ta20x5 | 1 | 16119 | 13512 | 17806 | 10676 | 10.47 | 20.99 |
| | 2 | 14693 | 15960 | 16096 | 12993 | 9.55 | 18.59 |
| | 3 | 16472 | 9329 | 16740 | 9060.8 | 1.63 | 2.87 |
| | 4 | 14436 | 17306 | 14926 | 15379 | 3.39 | 11.13 |
| | 5 | 15330 | 12978 | 15425 | 10391 | 0.62 | 19.93 |
| | 6 | 14729 | 6762.6 | 15172 | 6179 | 3.01 | 8.63 |
| | 7 | 15014 | 13057 | 16084 | 11335 | 7.13 | 13.19 |
| | 8 | 15947 | 14703 | 17472 | 10437 | 9.56 | 29.01 |
| | 9 | 14426 | 13147 | 15145 | 11280 | 4.98 | 14.20 |
| | 10 | 14502 | 13890 | 14650 | 12433 | 1.02 | 10.49 |
| | | | | | **Average** | **5.14** | **14.90** |
| Ta20x10 | 1 | 23154 | 78618 | 23585 | 74626 | 1.86 | 5.08 |
| | 2 | 24339 | 78259 | 25550 | 67715 | 4.98 | 13.47 |
| | 3 | 22184 | 81369 | 23066 | 69385 | 3.98 | 14.73 |
| | 4 | 20250 | 71046 | 20672 | 58327 | 2.08 | 17.90 |
| | 5 | 20778 | 55002 | 20811 | 51309 | 0.16 | 6.71 |
| | 6 | 21278 | 68738 | 21800 | 64052 | 2.45 | 6.82 |
| | 7 | 20399 | 65489 | 20707 | 60586 | 1.51 | 7.49 |
| | 8 | 22521 | 79259 | 23773 | 69351 | 5.56 | 12.50 |
| | 9 | 22520 | 60732 | 23157 | 55205 | 2.83 | 9.10 |
| | 10 | 23757 | 86507 | 24212 | 72129 | 1.92 | 16.62 |
| | | | | | **Average** | **2.73** | **11.04** |
| Ta20x20 | 1 | 34895 | 330990 | 37040 | 286371 | 6.15 | 13.48 |
| | 2 | 34809 | 283113 | 35473 | 262306 | 1.91 | 7.35 |
| | 3 | 35097 | 293074 | 36361 | 277031 | 3.60 | 5.47 |
| | 4 | 35080 | 302093 | 36693 | 277162 | 4.60 | 8.25 |
| | 5 | 36953 | 331550 | 37967 | 298250 | 2.74 | 10.04 |
| | 6 | 34440 | 278415 | 35591 | 250821 | 3.34 | 9.91 |
| | 7 | 34482 | 281575 | 35538 | 254530 | 3.06 | 9.60 |
| | 8 | 34699 | 279145 | 35876 | 266992 | 3.39 | 4.35 |
| | 9 | 35983 | 290376 | 36940 | 264060 | 2.66 | 9.06 |
| | 10 | 34988 | 293664 | 35839 | 278130 | 2.43 | 5.29 |
| | | | | | **Average** | **3.39** | **8.28** |
| | | | | | **Overall Average** | **3.75** | **11.41** |

*Contact: mallanjulap1@gmail.com*

**Table 12:** Taillard 50 jobs datasets with various machines configurations and their corresponding solutions with flowtime and energy consumption.

| Problem size (jobs x m/c) | Dataset | FT solution | | EC solution | | Variation in percentage | |
|---|---|---|---|---|---|---|---|
| | S.No | Min(FT1) | EC1 | FT2 | Min(EC2) | (FT2-FT1)/FT1 | (EC1-EC2)/EC1 |
| Ta50x5 | 1 | 81544 | 16369 | 87114 | 14265 | 6.83 | 12.85 |
| | 2 | 74635 | 17784 | 77001 | 14867 | 3.17 | 16.40 |
| | 3 | 81093 | 20669 | 86281 | 15444 | 6.40 | 25.28 |
| | 4 | 81955 | 30597 | 86000 | 25919 | 4.94 | 15.29 |
| | 5 | 83149 | 24104 | 86340 | 14394 | 3.84 | 40.28 |
| | 6 | 76056 | 25345 | 79975 | 20089 | 5.15 | 20.74 |
| | 7 | 78243 | 21958 | 84405 | 18287 | 7.88 | 16.72 |
| | 8 | 74320 | 23918 | 77087 | 22015 | 3.72 | 7.96 |
| | 9 | 72801 | 16967 | 75377 | 15903 | 3.54 | 6.27 |
| | 10 | 76690 | 26364 | 78417 | 20153 | 2.25 | 23.56 |
| | | | | | **Average** | **4.77** | **18.54** |
| Ta50x10 | 1 | 96104 | 93176 | 99736 | 79085 | 3.78 | 15.12 |
| | 2 | 97935 | 96770 | 102778 | 82257 | 4.95 | 15.00 |
| | 3 | 101663 | 102163 | 104568 | 83981 | 2.86 | 17.80 |
| | 4 | 103830 | 104332 | 107332 | 78297 | 3.37 | 24.95 |
| | 5 | 103398 | 110792 | 105032 | 86029 | 1.58 | 22.35 |
| | 6 | 104210 | 111253 | 109324 | 95000 | 4.91 | 14.61 |
| | 7 | 101890 | 109170 | 105992 | 99170 | 4.03 | 9.16 |
| | 8 | 100530 | 118491 | 104485 | 96949 | 3.93 | 18.18 |
| | 9 | 104991 | 105241 | 106841 | 95757 | 1.76 | 9.01 |
| | 10 | 100650 | 96926 | 104030 | 71084 | 3.36 | 26.66 |
| | | | | | **Average** | **3.45** | **17.28** |
| Ta50x20 | 1 | 136923 | 388912 | 144212 | 319463 | 5.32 | 17.86 |
| | 2 | 137263 | 366767 | 140821 | 350693 | 2.59 | 4.38 |
| | 3 | 138270 | 362951 | 138871 | 355423 | 0.43 | 2.07 |
| | 4 | 137152 | 359522 | 140653 | 330391 | 2.55 | 8.10 |
| | 5 | 138652 | 399572 | 140551 | 386719 | 1.37 | 3.22 |
| | 6 | 141812 | 437012 | 144167 | 392344 | 1.66 | 10.22 |
| | 7 | 140523 | 395085 | 145313 | 353653 | 3.41 | 10.49 |
| | 8 | 138832 | 383142 | 143342 | 341934 | 3.25 | 10.76 |
| | 9 | 136642 | 370164 | 136731 | 345692 | 0.07 | 6.61 |
| | 10 | 138280 | 346350 | 141150 | 313390 | 2.08 | 9.52 |
| | | | | | **Average** | **2.27** | **8.32** |
| | | | | | **Overall Average** | **3.50** | **14.71** |

*Contact: mallanjulap1@gmail.com*

**Table 13:** Taillard 100 jobs datasets with various machines configurations and their corresponding solutions with flowtime and energy consumption.

| Problem size (jobs x m/c) | Dataset S. No | Total flowtime Min(FT1) | EC1 | Energy consumption FT2 | Min(EC2) | Percentage variation (FT2-FT1)/FT1 | (EC1-EC2)/EC1 |
|---|---|---|---|---|---|---|---|
| Ta100x5 | 1 | 294163 | 20775 | 321121 | 15660 | 9.16 | 24.62 |
|  | 2 | 293241 | 17424 | 297443 | 14886 | 1.43 | 14.57 |
|  | 3 | 265332 | 30747 | 296271 | 19281 | 11.66 | 37.29 |
|  | 4 | 292191 | 28308 | 310471 | 19343 | 6.26 | 31.67 |
|  | 5 | 284990 | 37265 | 287182 | 20968 | 0.77 | 43.73 |
|  | 6 | 288142 | 39858 | 307142 | 27998 | 6.59 | 29.76 |
|  | 7 | 298520 | 19672 | 319892 | 14207 | 7.16 | 27.78 |
|  | 8 | 302031 | 24467 | 311090 | 15320 | 3.00 | 37.39 |
|  | 9 | 296412 | 28817 | 301880 | 23166 | 1.84 | 19.61 |
|  | 10 | 297390 | 24759 | 299880 | 23514 | 0.84 | 5.03 |
|  |  |  |  |  | **Average** | **4.87** | **27.14** |
| Ta100x10 | 1 | 366850 | 147591 | 373942 | 124152 | 1.93 | 15.88 |
|  | 2 | 338132 | 149912 | 356761 | 110623 | 5.51 | 26.21 |
|  | 3 | 350912 | 149060 | 352316 | 129112 | 0.40 | 13.38 |
|  | 4 | 364914 | 125013 | 377713 | 106384 | 3.51 | 14.90 |
|  | 5 | 350253 | 124712 | 354052 | 115890 | 1.08 | 7.07 |
|  | 6 | 331572 | 129390 | 343780 | 112222 | 3.68 | 13.27 |
|  | 7 | 343850 | 160812 | 345363 | 148790 | 0.44 | 7.48 |
|  | 8 | 351010 | 128542 | 353181 | 119854 | 0.62 | 6.76 |
|  | 9 | 342471 | 149712 | 358572 | 139732 | 4.70 | 6.67 |
|  | 10 | 355213 | 141423 | 361950 | 111160 | 1.90 | 21.40 |
|  |  |  |  |  | **Average** | **2.38** | **13.30** |
| Ta100x20 | 1 | 441050 | 512132 | 455460 | 475290 | 3.27 | 7.19 |
|  | 2 | 443370 | 473769 | 454490 | 447280 | 2.51 | 5.59 |
|  | 3 | 434340 | 470960 | 440613 | 437340 | 1.44 | 7.14 |
|  | 4 | 436860 | 516123 | 450238 | 452848 | 3.06 | 12.26 |
|  | 5 | 439512 | 508043 | 441652 | 493190 | 0.49 | 2.92 |
|  | 6 | 437669 | 480480 | 439183 | 466823 | 0.35 | 2.84 |
|  | 7 | 435332 | 488312 | 440310 | 414230 | 1.14 | 15.17 |
|  | 8 | 447230 | 544020 | 463140 | 436456 | 3.56 | 19.77 |
|  | 9 | 443612 | 574672 | 453521 | 474670 | 2.23 | 17.40 |
|  | 10 | 448923 | 573762 | 458221 | 504452 | 2.07 | 12.08 |
|  |  |  |  |  | **Average** | **2.01** | **10.24** |
|  |  |  |  |  | **Overall Average** | **3.09** | **16.89** |

*Contact: mallanjulap1@gmail.com*

## 5. Conclusions and future scope

In this paper, we investigated permutation flowshop scheduling problem that involves both sustainability related criterion (i.e., the energy consumption) and the productivity related criterion (i.e., the flowtime). Integrating the eco-friendly contemplation into scheduling verdict in such a way can directly minimize the energy consumption and flowtime in the manufacturing enterprises.

The discovered fronts generated by this algorithm can provide useful information for decision makers to achieve good trade-off between the flowtime and the energy consumption in production scheduling. Computational experiments are conducted on the benchmark problems of Taillard sets, to reveal the effectiveness of the above-mentioned procedure and algorithm. The results demonstrated the importance of parameter tuning in finding the best genetic operator values in evaluating the performance of the proposed algorithm with the literature article. The results also reveal that the best parameters selection of the NSGA-II algorithm outperforms the literature papers in both solution quality and diversity. Finally, the novel contribution of this research can be given in three folds. First, we constructed a mathematical model that reduces the energy consumption of machine tools and simultaneously minimizes the flowtime is proposed. Second, we attempted to conduct the exhaustive experiments in order to improve the proposed NSGA-II algorithm in parameter tuning to find the best optimal solution when compared to the literature solutions. Third, we have solved each Taillard dataset 10 times to come-up with benchmark results using fine-tuned meta-heuristic algorithm. These results can act as meta-heuristic benchmark results for Taillard flowshop problem by considering energy as a parameter. Future researchers can consider these results as baseline results when they use heuristics or meta-heuristics.

*Contact: mallanjulap1@gmail.com*

# References

Afshin Mansouri, S., & Aktas, E. (2016). Minimizing energy consumption and makespan in a two-machine flowshop scheduling problem. *Journal of the Operational Research Society*, *67*(11), 1382–1394.

Ahilan, C., Kumanan, S., Sivakumaran, N., & Edwin Raja Dhas, J. (2013). Modeling and prediction of machining quality in CNC turning process using intelligent hybrid decision making tools. *Applied Soft Computing Journal*, *13*(3). https://doi.org/10.1016/j.asoc.2012.03.071

Akbar, M., & Irohara, T. (2018). Scheduling for sustainable manufacturing: A review. *Journal of Cleaner Production*.

Amiri, M. F., & Behnamian, J. (2020). Multi-objective green flowshop scheduling problem under uncertainty: Estimation of distribution algorithm. *Journal of Cleaner Production*, *251*, 119734.

Artigues, C., Lopez, P., & Hait, A. (2013). The energy scheduling problem: Industrial case-study and constraint propagation techniques. *International Journal of Production Economics*, *143*(1), 13–23.

Che, A., Lv, K., Levner, E., & Kats, V. (2015). Energy consumption minimization for single machine scheduling with bounded maximum tardiness. *ICNSC 2015 - 2015 IEEE 12th International Conference on Networking, Sensing and Control*. https://doi.org/10.1109/ICNSC.2015.7116025

Che, A., Wu, X., Peng, J., & Yan, P. (2017). Energy-efficient bi-objective single-machine scheduling with power-down mechanism. *Computers & Operations Research*, *85*, 172–183.

Dai, M., Tang, D., Giret, A., Salido, M. A., & Li, W. D. (2013). Energy-efficient scheduling for a flexible flow shop using an improved genetic-simulated annealing algorithm. *Robotics and Computer-Integrated Manufacturing*, *29*(5), 418–429.

Deb, K., Pratap, A., Agarwal, S., & Meyarivan, T. (2002). A fast and elitist multiobjective genetic algorithm: NSGA-II. *IEEE Transactions on Evolutionary Computation*, *6*(2), 182–197.

Diaz, N., Redelsheimer, E., & Dornfeld, D. (2011). Energy consumption characterization and reduction strategies for milling machine tool use. *Glocalized Solutions for Sustainability in Manufacturing*, 263–267.

Ding, J.-Y., Song, S., & Wu, C. (2016). Carbon-efficient scheduling of flow shops by multi-objective optimization. *European Journal of Operational Research*, *248*(3), 758–771.

Duflou, J. R., Sutherland, J. W., Dornfeld, D., Herrmann, C., Jeswiet, J., Kara, S., Hauschild, M., & Kellens, K. (2012). Towards energy and resource efficient manufacturing: A processes and systems approach. *CIRP Annals*, *61*(2), 587–609.

Fang, K., Uhan, N. A., Zhao, F., & Sutherland, J. W. (2013). Flow shop scheduling with peak power consumption constraints. *Annals of Operations Research*, *206*(1), 115–145.

Fang, K., Uhan, N., Zhao, F., & Sutherland, J. W. (2011a). A new approach to scheduling in manufacturing for power consumption and carbon footprint reduction. *Journal of Manufacturing Systems*, *30*(4), 234–240.

Fang, K., Uhan, N., Zhao, F., & Sutherland, J. W. (2011b). A new approach to scheduling in manufacturing for power consumption and carbon footprint reduction. *Journal of Manufacturing Systems*, *30*(4), 234–240.

Gahm, C., Denz, F., Dirr, M., & Tuma, A. (2016). Energy-efficient scheduling in manufacturing companies: A review and research framework. *European Journal of Operational Research*, *248*(3), 744–757.

Gharbi, A., Ladhari, T., Msakni, M. K., & Serairi, M. (2013a). The two-machine flowshop scheduling problem with sequence-independent setup times: New lower bounding strategies. *European Journal of Operational Research*, *231*(1), 69–78.

Gharbi, A., Ladhari, T., Msakni, M. K., & Serairi, M. (2013b). The two-machine flowshop scheduling problem with sequence-independent setup times: New lower bounding strategies. *European Journal of Operational Research*, *231*(1), 69–78.

Ghazouani, A., Jebli, M. ben, & Shahzad, U. (2021). Impacts of environmental taxes and technologies on greenhouse gas emissions: contextual evidence from leading emitter European countries. *Environmental Science and Pollution Research*, *28*(18), 22758–22767.

Gunasekaran, A., & Ngai, E. W. T. (2012). The future of operations management: an outlook and analysis. *International Journal of Production Economics*, *135*(2), 687–701.

Haapala, K. R., Zhao, F., Camelio, J., Sutherland, J. W., Skerlos, S. J., Dornfeld, D. A., Jawahir, I. S., Clarens, A. F., & Rickli, J. L. (2013). A review of engineering research in sustainable manufacturing. *Journal of Manufacturing Science and Engineering*, *135*(4).

He, Y., Liu, F., Cao, H., & Li, C. (2005). A bi-objective model for job-shop scheduling problem to minimize both energy consumption and makespan. *Journal of Central South University of Technology*, *12*(2), 167–171.

Kizilay, D., Tasgetiren, M. F., Pan, Q.-K., & Süer, G. (2019). An ensemble of meta-heuristics for the energy-efficient blocking flowshop scheduling problem. *Procedia Manufacturing*, *39*, 1177–1184.

Kleindorfer, P. R., Singhal, K., & van Wassenhove, L. N. (2005). Sustainable operations management. *Production and Operations Management*, *14*(4), 482–492.

Lei, H., Wang, R., Zhang, T., Liu, Y., & Zha, Y. (2016). A multi-objective co-evolutionary algorithm for energy-efficient scheduling on a green data center. *Computers & Operations Research*, *75*, 103–117.

Li, J., Sang, H., Han, Y., Wang, C., & Gao, K. (2018). Efficient multi-objective optimization algorithm for hybrid flow shop scheduling problems with setup energy consumptions. *Journal of Cleaner Production*, *181*, 584–598.

Li, S., Liu, F., & Zhou, X. (2018). Multi-objective energy-saving scheduling for a permutation flow line. *Proceedings of the Institution of Mechanical Engineers, Part B: Journal of Engineering Manufacture*, *232*(5), 879–888.

Li, X., Xing, K., Wu, Y., Wang, X., & Luo, J. (2017). Total energy consumption optimization via genetic algorithm in flexible manufacturing systems. *Computers & Industrial Engineering*, *104*, 188–200.

*Contact: mallanjulap1@gmail.com*

Liang, P., Yang, H., Liu, G., & Guo, J. (2015). An ant optimization model for unrelated parallel machine scheduling with energy consumption and total tardiness. *Mathematical Problems in Engineering*, *2015*.

Lin, S.-W., Ying, K.-C., & Lee, Z.-J. (2009). Metaheuristics for scheduling a non-permutation flowline manufacturing cell with sequence dependent family setup times. *Computers & Operations Research*, *36*(4), 1110–1121.

Lin, W., Yu, D. Y., Zhang, C., Liu, X., Zhang, S., Tian, Y., Liu, S., & Xie, Z. (2015). A multi-objective teaching–learning-based optimization algorithm to scheduling in turning processes for minimizing makespan and carbon footprint. *Journal of Cleaner Production*, *101*, 337–347.

Liu, G.-S., Zhang, B.-X., Yang, H.-D., Chen, X., & Huang, G. Q. (2013). A branch-and-bound algorithm for minimizing the energy consumption in the PFS problem. *Mathematical Problems in Engineering*, *2013*.

Liu, W., Jin, Y., & Price, M. (2017). A new improved NEH heuristic for permutation flowshop scheduling problems. *International Journal of Production Economics*, *193*, 21–30.

Mansouri, S. A., Aktas, E., & Besikci, U. (2016). Green scheduling of a two-machine flowshop: Trade-off between makespan and energy consumption. *European Journal of Operational Research*, *248*(3), 772–788.

Mashaei, M., & Lennartson, B. (2012). Energy reduction in a pallet-constrained flow shop through on–off control of idle machines. *IEEE Transactions on Automation Science and Engineering*, *10*(1), 45–56.

Masmoudi, O., Yalaoui, A., Ouazene, Y., & Chehade, H. (2016). Multi-item capacitated lot-sizing problem in a flow-shop system with energy consideration. *IFAC-PapersOnLine*, *49*(12), 301–306.

Moon, J.-Y., Shin, K., & Park, J. (2013). Optimization of production scheduling with time-dependent and machine-dependent electricity cost for industrial energy efficiency. *The International Journal of Advanced Manufacturing Technology*, *68*(1), 523–535.

Mourtzis, D., Vlachou, E., Milas, N., & Dimitrakopoulos, G. (2016). Energy consumption estimation for machining processes based on real-time shop floor monitoring via wireless sensor networks. *Procedia CIRP*, *57*, 637–642.

Mouzon, G., & Yildirim, M. B. (2008). A framework to minimise total energy consumption and total tardiness on a single machine. *International Journal of Sustainable Engineering*, *1*(2), 105–116.

Neto, R. F. T., & Godinho Filho, M. (2011). An ant colony optimization approach to a permutational flowshop scheduling problem with outsourcing allowed. *Computers & Operations Research*, *38*(9), 1286–1293.

Ngai, E. W. T., To, C. K. M., Ching, V. S. M., Chan, L. K., Lee, M. C. M., Choi, Y. S., & Chai, P. Y. F. (2012). Development of the conceptual model of energy and utility management in textile processing: A soft systems approach. *International Journal of Production Economics*, *135*(2), 607–617.

Nilakantan, J. M., Li, Z., Tang, Q., & Nielsen, P. (2017). Multi-objective co-operative co-evolutionary algorithm for minimizing carbon footprint and maximizing line efficiency in robotic assembly line systems. *Journal of Cleaner Production*, *156*, 124–136.

Rager, M., Gahm, C., & Denz, F. (2015). Energy-oriented scheduling based on evolutionary algorithms. *Computers & Operations Research*, *54*, 218–231.

Saddikuti, V., & Pesaru, V. (2019). NSGA Based Algorithm for Energy Efficient Scheduling in Cellular Manufacturing. *Procedia Manufacturing*, *39*, 1002–1009. https://doi.org/https://doi.org/10.1016/j.promfg.2020.01.379

Salgado, F., & Pedrero, P. (2008). Short-term operation planning on cogeneration systems: A survey. *Electric Power Systems Research*, *78*(5), 835–848.

Salido, M. A., Escamilla, J., Barber, F., Giret, A., Tang, D., & Dai, M. (2013). Energy-aware parameters in job-shop scheduling problems. *GREEN-COPLAS 2013: IJCAI 2013 Workshop on Constraint Reasoning, Planning and Scheduling Problems for a Sustainable Future*, 44–53.

Schulz, S., Schönheit, M., & Neufeld, J. S. (2022). Multi-objective carbon-efficient scheduling in distributed permutation flow shops under consideration of transportation efforts. *Journal of Cleaner Production*, 132551.

Shrouf, F., Ordieres-Meré, J., García-Sánchez, A., & Ortega-Mier, M. (2014). Optimizing the production scheduling of a single machine to minimize total energy consumption costs. *Journal of Cleaner Production*, *67*, 197–207.

Song, W. J., Zhang, C. Y., Lin, W. W., & Shao, X. Y. (2014). Flexible job-shop scheduling problem with maintenance activities considering energy consumption. *Applied Mechanics and Materials*, *521*, 707–713.

Taillard, E. (1993). Benchmarks for basic scheduling problems. *European Journal of Operational Research*, *64*(2), 278–285.

Tang, D., Dai, M., Salido, M. A., & Giret, A. (2016). Energy-efficient dynamic scheduling for a flexible flow shop using an improved particle swarm optimization. *Computers in Industry*, *81*, 82–95.

Wang, X., & Tang, L. (2017). A machine-learning based memetic algorithm for the multi-objective permutation flowshop scheduling problem. *Computers & Operations Research*, *79*, 60–77.

Wu, T.-Y., Wu, I.-C., & Liang, C.-C. (2013). Multi-objective flexible job shop scheduling problem based on Monte-Carlo tree search. *2013 Conference on Technologies and Applications of Artificial Intelligence*, 73–78.

Yang, X., Zeng, Z., Wang, R., & Sun, X. (2016a). Bi-objective flexible job-shop scheduling problem considering energy consumption under stochastic processing times. *PloS One*, *11*(12), e0167427.

Yang, X., Zeng, Z., Wang, R., & Sun, X. (2016b). Bi-objective flexible job-shop scheduling problem considering energy consumption under stochastic processing times. *PloS One*, *11*(12), e0167427.

Yüksel, D., Taşgetiren, M. F., Kandiller, L., & Gao, L. (2020a). An energy-efficient bi-objective no-wait permutation flowshop scheduling problem to minimize total tardiness and total energy consumption. *Computers & Industrial Engineering*, *145*, 106431.

Yüksel, D., Taşgetiren, M. F., Kandiller, L., & Gao, L. (2020b). An energy-efficient bi-objective no-wait permutation flowshop scheduling problem to minimize total tardiness and total energy consumption. *Computers & Industrial Engineering*, *145*, 106431.

Zanoni, S., Bettoni, L., & Glock, C. H. (2014). Energy implications in a two-stage production system with controllable production rates. *International Journal of Production Economics*, *149*, 164–171.

Zhang, C., Gu, P., & Jiang, P. (2015). Low-carbon scheduling and estimating for a flexible job shop based on carbon footprint and carbon efficiency of multi-job processing. *Proceedings of the Institution of Mechanical Engineers, Part B: Journal of Engineering Manufacture*, *229*(2), 328–342.

Zhang, H., Zhao, F., Fang, K., & Sutherland, J. W. (2014). Energy-conscious flow shop scheduling under time-of-use electricity tariffs. *CIRP Annals*, *63*(1), 37–40.

Zhang, R., & Chiong, R. (2016). Solving the energy-efficient job shop scheduling problem: A multi-objective genetic algorithm with enhanced local search for minimizing the total weighted tardiness and total energy consumption. *Journal of Cleaner Production*, *112*, 3361–3375.

*Contact: mallanjulap1@gmail.com*